\documentclass{ifacconf}

\usepackage{graphicx}      
\usepackage{natbib}        
\usepackage{upgreek}
\usepackage{amsfonts}
\usepackage{amsmath}
\usepackage[many]{tcolorbox}
\usepackage{graphicx} 
\newtheorem{lemma}[thm]{Lemma}
\newtheorem{remark}{Remark}
\newtheorem{definition}[thm]{Definition}

\begin{document}
\begin{frontmatter}

\title{Characterization of optimization problems that are solvable iteratively with linear convergence } 


\author[First]{Foivos Alimisis}

\address[First]{University of Geneva, 
   Switzerland}

\begin{abstract}                
In this work, we state a general conjecture on the solvability of optimization problems
via algorithms with linear convergence guarantees. We make a first step towards examining its
correctness by fully characterizing the problems that are solvable via Riemannian gradient
descent with linear convergence.
\end{abstract}

\begin{keyword}
Optimization; Riemannian optimization; Gradient descent; Riemannian gradient descent; (Geodesic) weak-strong-convexity.
\end{keyword}

\end{frontmatter}

\section{Introduction}
\vspace{-2mm}
In recent years, advances in the theory of optimization on manifolds [\cite{absilOptimizationAlgorithmsMatrix2008,boumal2023introduction}] have highlighted the importance of Riemannian algorithms in solving problems that are generally intractable if posed in the Euclidean setting. Such a problem is for instance operator scaling [\cite{gurvits2004classical}]. This problem is not convex if posed in Euclidean space, but it is ``geodesically" convex if posed in the space of symmetric and positive definite matrices with a carefully chosen Riemannian metric [\cite{allen2018operator}]. This allows the design of Riemannian algorithms that extend algorithms from the realm of convex optimization and can solve such a problem with theoretical guarantees [\cite{zhang2016first,zhang2018towards,pmlr-v162-kim22k}]. 

More recently \cite{alimisis2022geodesic} showed that the symmetric eigenvalue problem (which was long known to be solvable with linear convergence via subspace iteration) is geodesically weak-strongly-convex on the Grassmann manifold. Weak-strong-convexity is a property weaker than strong convexity that still allows to analyse the convergence of many algorithms. More precisely, it is enough to guarantee linear convergence of gradient descent under a carefully chosen step size, as it is showed in [\cite{alimisis2022geodesic}] influenced by the results in [\cite{bu2020note}]. For a nice presentation of weak-strong-convexity and its relationship to other important convexity notions, the reader is referred to Appendix A in [\cite{karimi2016linear}]. A slightly different but close in the spirit study can be found in [\cite{necoara2019linear,zhang2019acceleration}]. 
We note that there are some small differences in the definition of weak-strong-convexity throughout the literature, but the spirit is always the same. 

The previously discussed progress is reflected in the following ambitious conjecture:
\begin{tcolorbox}
If the solution of an optimization problem can be approximated iteratively by an algorithm with linear convergence rate, then the objective cost function must satisfy some convexity-like property with respect to some Riemannian metric.
\end{tcolorbox}
As it is outlined in the rest of this work, there are reasons to believe that this convexity-like property must be weak-strong-convexity. 
Our main contribution is to prove that any smooth optimization problem that can be solved via gradient descent with linear convergence must necessarily be weak-strongly-convex (Theorem \ref{thm:main}). This result can be extended to an optimization problem in a Riemannian manifold of bounded sectional curvatures solved via Riemannian gradient descent (Theorem \ref{thm:main_riemannian}). This implies that the previous conjecture is true if the algorithm used to solve the objective optimization problem can be written as Riemannian gradient descent with respect to some Riemannian metric. This is probably not possible for any algorithm, but the free choice of any Riemannian metric provides a lot of options. For instance, standard Euclidean preconditioned gradient descent can be written as Riemannian gradient descent with respect to a flat ($0$ curvature) Riemannian metric. We treat this case in Appendix C.

\section{Important notions and background}
For the rest, $f$ will be a $C^1$ real-valued function, which will be our optimization objective. For now, it will be defined in an \emph{open and convex} subset $E$ of the Euclidean space. 

We proceed with the definition of weak-strong-convexity. This notion is equivalent with weak-quasi-convexity [\cite{hardt2018gradient,guminov2017accelerated}] and the more well-known quadratic growth condition [\cite{drusvyatskiy2018error}] holding simultaneously.
\begin{definition}
\label{def:weak-strong-conv}
   A (lower bounded) function $f: E  \rightarrow \mathbb{R}$, where $E \subseteq \mathbb{R}^n$ is open and convex, is called weak-strongly-convex, if for \emph{any} global minimizer $x^*$ and any point $x \in E$, it holds
\begin{equation}
\label{eq:weak-quasi-conv}
    f(x)-f(x^*) \leq \frac{1}{a} \langle \nabla f(x), x-x^* \rangle - \frac{\mu}{2} \|x-x^*\|^2,
\end{equation}
where $a,\mu>0$ are constants. 
\end{definition}


 Inequality (\ref{eq:weak-quasi-conv}) implies that there is only one global minimizer $x^* \in E$, because if there were global minimizers $x^*_1 \neq x^*_2$, then
\begin{equation*}
    \underbrace{f(x^*_1)-f(x^*_2)}_{=0} \leq \frac{1}{a} \underbrace{\langle \nabla f(x^*_1), x^*_1-x^*_2 \rangle}_{=0} - \underbrace{\frac{\mu}{2} \|x^*_1-x^*_2\|^2}_{>0},
\end{equation*}
which cannot hold. 

From now on, we denote the unique global minimizer of $f$ by $x^*$. 

Also, there are no other critical points:
indeed, if $y \in E$ is a critical point of $f$, then 
\begin{equation*}
    \underbrace{f(y)-f(x^*)}_{\geq 0} \leq \frac{1}{a} \underbrace{\langle \nabla f(y), y-x^* \rangle}_{=0} - \frac{\mu}{2} \underbrace{\|y-x^*\|^2}_{\geq 0},
\end{equation*}
thus $y=x^*$.

We now pass to the definition of smoothness, which is the most classic assumption imposed to optimization problems (and satisfied in plenty of practical situations):
\begin{definition}
\label{def:smoothness}
    $f:E \rightarrow \mathbb{R}$ is called $\gamma$-smooth if its gradient is $\gamma$-Lipschitz continuous, i.e.
\begin{equation*}
    \| \nabla f(x)-\nabla f(y) \| \leq \gamma \|x-y\|,
\end{equation*}
for any $x,y \in E$.
\end{definition}

It is known [\cite{bu2020note}, Lemma 4.2]  that if a function is weak-strongly-convex and smooth in $\mathbb{R}^n$, then gradient descent implemented with a suitable step size provides a contraction rate at each step. 

We denote one step of gradient descent starting from $x \in \mathbb{R}^n$ as
\begin{equation}
\label{eq:GD}
    \tilde x := x-\eta \nabla f(x).
\end{equation}
Here $\eta$ is the step size.

We provide Lemma 4.2 of \cite{bu2020note} in the following proposition:
\begin{prop}\label{prop:convergence}
Consider the problem of minimizing a $\gamma$-smooth and weak-strongly-convex function $f:\mathbb{R}^n \rightarrow \mathbb{R}$ with parameters $a$ and $\mu$ as in inequality (\ref{eq:weak-quasi-conv}). If $\tilde x$ is produced by equation (\ref{eq:GD}) with $\eta \leq a/\gamma$, then we have
\begin{equation*}
    \|\tilde x-x^* \|^2 \leq (1-a \mu \eta) \|x-x^*\|^2.
\end{equation*}
\end{prop}

\begin{remark}
   In fact, we do not need $f$ to be $\gamma$-smooth but only to saisfy the weaker property
\begin{equation}
\label{eq:annonying_weaker_smooth}
    f(x)-f(x^*) \geq \frac{1}{2\gamma}\|\nabla f(x)\|^2.
\end{equation}
Proposition \ref{prop:convergence} states that, in contrast to the popular belief that strong convexity is crucial to achieve strong convergence guarantees for gradient descent, a weaker star-type property like weak-strong-convexity (holds only around the optimum $x^*$) is all we need.  
\end{remark} 

\begin{remark}
It is interesting to note that Proposition \ref{prop:convergence} holds if $\mathbb{R}^n$ is substituted by any open and convex subset of it, if the $\gamma$-smoothness assumption is substituted by inequality (\ref{eq:annonying_weaker_smooth}). Also note that when $E=\mathbb{R}^n$, inequality (\ref{eq:annonying_weaker_smooth}) is implied by $\gamma$-smoothness, but when $E$ is a strict subset of $\mathbb{R}^n$ this is not true in general.
\end{remark}

\section{Main result}
\label{euclidean_result}
As discussed in the introduction, the main objective of this paper is to prove a version of the inverse of Proposition \ref{prop:convergence}. We state our result for an optimization problem defined in an open and convex subset $E \subseteq \mathbb{R}^n$, which could be the whole $\mathbb{R}^n$.

Before stating our main result (Theorem \ref{thm:main}), we need a property that is implied by $\gamma$-smoothness:
\begin{lemma}
\label{le:weaker_smoothness}
    If $E \subseteq \mathbb{R}^n$ is open and convex and $f:E \rightarrow \mathbb{R}$ is $\gamma$-smooth, then
    \begin{equation*}
        \| x-x^*\|^2 \geq \frac{2}{\gamma} (f(x)-f(x^*)).
    \end{equation*}
\end{lemma}
\begin{pf}
    The proof is standard, the reader can refer to the nice presentation of Lemma 4 in [\cite{zhou2018fenchel}] for details (take $y=x^*$ in point [2]). \qed
\end{pf}

\begin{tcolorbox}
\begin{thm}
\label{thm:main}
    Let $E \subseteq \mathbb{R}^n$ be open and convex. Consider the optimization problem
    \begin{equation*}
        \min_{x \in E} f(x),
    \end{equation*}
where $f\colon E \to \mathbb{R}$ is $C^1$, $\gamma$-smooth and has a unique global minimizer $x^*$. \newline
If a step of gradient descent  with step-size $\eta$ starting from any point $x \in E$  satisfies
\begin{equation*}
    \| \tilde x-x^* \|^2 \leq (1-c) \|x-x^*\|^2,
\end{equation*}
for a constant $0<c \leq 1$, then $f$ is weak-strongly-convex. The parameters of inequality \eqref{eq:weak-quasi-conv} satisfy
\begin{equation*}
    a:=\frac{1}{2 \gamma \eta} \frac{c}{1-\sqrt{c}/2} = \Omega \left(\frac{c}{\gamma \eta} \right), \hspace{2mm}\mu := \frac{\gamma}{2} = \Omega(\gamma).
\end{equation*}

\end{thm}
\end{tcolorbox}
\begin{pf}
Let $x \in E$ and $\tilde x$ the result of one iteration of gradient descent (\ref{eq:GD}).

    We first rewrite the term $\| \tilde x-x^* \|^2$:
    \begin{align*}
        \| \tilde x-x^* \|^2 &=\| x-\eta\nabla f(x)-x^* \|^2 \\&= \|x-x^*\|^2 -2\eta \langle \nabla f(x),x-x^* \rangle + \eta^2 \|\nabla f(x)\|^2.
    \end{align*}
This equality together with the contraction assumption gives
\begin{equation*}
    \|x-x^*\|^2 -2\eta \langle \nabla f(x),x-x^* \rangle + \eta^2 \|\nabla f(x)\|^2 \leq (1-c) \|x-x^*\|^2,
\end{equation*}
which can be rewritten as
\begin{equation}
\label{eq:2}
    2\eta \langle \nabla f(x),x-x^* \rangle \geq c \|x-x^*\|^2 + \eta^2 \|\nabla f(x)\|^2.
\end{equation}
Next, we use the inequality
\begin{equation*}
    \langle \alpha,\beta \rangle \leq \frac{\rho}{2} \|\alpha\|^2+\frac{1}{2 \rho} \| \beta \|^2
\end{equation*}
for any $\alpha,\beta \in \mathbb{R}^n$ and any $\rho>0$ to obtain
\begin{equation*}
    \frac{\rho}{2} \| \nabla f(x)\|^2 \geq \langle \nabla f(x),x-x^* \rangle-\frac{1}{2 \rho} \|x-x^*\|^2.
\end{equation*}
Multiplying both sides by $\frac{2 \eta^2}{\rho}$, we get
\begin{equation*}
    \eta^2 \| \nabla f(x)\|^2 \geq \frac{2 \eta^2}{\rho}\langle \nabla f(x),x-x^* \rangle-\frac{\eta^2}{\rho^2} \|x-x^*\|^2.
\end{equation*}
Using equation \eqref{eq:2}, we get
\begin{align*}
2\eta \langle \nabla f(x),x-x^* \rangle \geq c \|x-x^*\|^2 & + \frac{2 \eta^2}{\rho}\langle \nabla f(x),x-x^* \rangle \\&-\frac{\eta^2}{\rho^2} \|x-x^*\|^2,
\end{align*}
or equivalently
\begin{equation*}
\left(2\eta - \frac{2 \eta^2}{\rho} \right) \langle \nabla f(x),x-x^* \rangle \geq \left( c-\frac{\eta^2}{\rho^2} \right) \|x-x^*\|^2.
\end{equation*}
Since the last inequality holds for any $\rho>0$, we can choose $\rho= 2 \frac{\eta}{\sqrt{c}}$. Then it becomes
\begin{align*}
    2 \eta \left( 1-\frac{\sqrt{c}}{2} \right) \langle \nabla f(x),x-x^* \rangle &\geq \frac{3 c}{4} \|x-x^*\|^2 \\&= \frac{c}{4} \|x-x^*\|^2 + \frac{c}{2} \|x-x^*\|^2.
\end{align*}
By Lemma \ref{le:weaker_smoothness}, we have
\begin{equation*}
    \|x-x^*\|^2 \geq \frac{2}{\gamma} (f(x)-f(x^*)),
\end{equation*}
and using that to bound the last term of the previous inequality, we have
\begin{equation*}
    2 \eta \left( 1-\frac{\sqrt{c}}{2} \right) \langle \nabla f(x),x-x^* \rangle \geq \frac{c}{4} \|x-x^*\|^2 + \frac{c}{\gamma} (f(x)-f(x^*)).
\end{equation*}
Rearranging, we get
\begin{equation*}
    f(x)-f(x^*) \leq \underbrace{2 \gamma \eta \frac{ 1-\frac{\sqrt{c}}{2} }{c 
 }}_{\frac{1}{a}} \langle \nabla f(x),x-x^* \rangle - \underbrace{\frac{\gamma}{4}}_{\frac{\mu}{2}} \|x-x^*\|^2.
\end{equation*}
The last inequality provides the desired result. \qed
\end{pf}

\begin{remark}
In order for the result of Theorem \ref{thm:main} to make sense, we need to have $c \geq a \mu \eta$, otherwise we would be able to apply Proposition \ref{prop:convergence} and improve the assumed convergence rate for free. Ideally, we would need $c=a \mu \eta$, which means that the result is optimal. In our case we have
\begin{equation*}
    a \mu \eta = \frac{1}{4} \frac{c}{1-\frac{\sqrt{c}}{2}},
\end{equation*}
which is between $\frac{c}{4}$ and $\frac{c}{2}$. This means that our result is reasonable, but slightly suboptimal.
\end{remark} 

Theorem \ref{thm:main} states that a smooth function for which gradient descent converges in a steady contraction rate with respect to the distance of the iterates to the optimum (and the iterates remain inside $E$) is necessarily weak-strongly-convex. This means that weak-strong-convexity is in some sense necessary and sufficient if one wishes to obtain a linear convergence rate for gradient descent.

\section{Extension of Theorem \ref{thm:main} to manifold optimization}
\label{sec:riemannian}

In this section, we extend our main result to Riemannian gradient descent in a complete Riemannian manifold $M$ of sectional curvatures bounded from above. Due to lack of space, we expose some basic notions from Riemannian geometry in Appendix A. The reader can refer to it for basic definitions and notations that is used in the rest of this work.

We begin by providing the definition of geodesic weak-strong-convexity as it appears in [\cite{alimisis2022geodesic}]. This is an analogue of Definition \ref{def:weak-strong-conv}.

\begin{definition}
\label{def:weak-strong-conv_riem}
    Let $E \subseteq M$ open and geodesically convex. A (lower bounded) function $f: E \rightarrow \mathbb{R}$ is called geodesically weak-strongly-convex, if for \emph{any} global minimizer $x^*$ and any point $x \in E$, it holds
\begin{equation}
\label{eq:weak-quasi-conv_riem}
    f(x)-f(x^*) \leq \frac{1}{a} \langle \textnormal{grad}f(x), -\log_x(x^*) \rangle - \frac{\mu}{2} \textnormal{dist}^2(x,x^*),
\end{equation}
where $a,\mu>0$ are constants.
\end{definition}

As in the Euclidean case, it is easy to show that a geodesically weak-strongly-convex function has only one critical point and this is the global minimizer. We denote this global minimizer again by $x^*$.

We now discuss the notion of geodesic smoothness:
\begin{definition}
\label{def:geod_smoothness}
A function $f:E \rightarrow \mathbb{R}$, $E \subseteq M$ open and geodesically convex, is called geodesically $\gamma$-smooth if for any $x,y \in E$ it holds
\begin{equation}
\label{eq:geod-smoothness}
    \| \textnormal{\textnormal{gradf}}(x) -\Gamma_y^x {\textnormal{\textnormal{gradf}}(y)} \| \leq \gamma \textnormal{dist}(x,y),
\end{equation}
where $\Gamma_y^x$ is the parallel transport from $y$ to $x$ along the geodesic that connects them.
\end{definition}
This definition is a natural Riemannian extension of Definition \ref{def:smoothness} that the reader can find in the vast majority of works on Riemannian algorithms, see for instance [\cite{zhang2016first}, Definition 4] or [\cite{alimisis2020continuous}, Definition 9]. It is particularly natural since, in the case that $f$ is twice differentiable, $\gamma$ is an upper bound of the Riemannian analogue of the Hessian.

When $E=M$, these properties are enough to guarantee linear convergence for a version of Riemannian gradient descent, which reads as
\begin{equation}
\label{eq:riemannian_gd}
    \tilde x = \exp_x(-\eta \textnormal{grad}f(x))
\end{equation}
with $\eta$ being a step-size. This is evident by Lemma 11 in [\cite{alimisis2022geodesic}], in the case that the sectional curvatures of the manifold are all nonnegative (Lemma 11 in [\cite{alimisis2022geodesic}] is applied to the specific case of the Grassmann manifold, but it can be easily adapted to any manifold of nonnegative sectional curvatures). Moreover, it is easy to prove a similar result in the case that $M$ can have negative curvatures, by utilizing the (now classic) Lemma 5 in [\cite{zhang2016first}] and implementing Riemannian gradient descent with stepsize $\eta \leq \frac{a}{\zeta \gamma}$. Here $\zeta$ is the curvature-dependent quantity appearing in Lemma 5 of [\cite{zhang2016first}].

For completeness, we present these results unified and in the slightly more general case that $f$ is defined in an open and geodesically convex subset $E$ of $M$. First we need a small reformulation of geodesic smoothness:

Assumption:
    $f:E (\subseteq M) \rightarrow \mathbb{R}$ is assumed to satisfy the inequality
    \begin{equation}
    \label{eq:annonying_weaker_smooth_riem}
        f(x)-f(x^*) \geq \frac{1}{2 \gamma} \| \textnormal{grad}f(x) \|^2.
    \end{equation}

Note that if $E=M$, then inequality (\ref{eq:annonying_weaker_smooth_riem}) is implied by geodesic $\gamma$-smoothness. If $E$ is a strict subset of $M$ this may be incorrect.




\begin{prop}
\label{prop:convergence_riem}
Let $E$ be an open and geodesically convex subset of $M$, with $M$ being a Riemannian manifold of sectional curvatures bounded from below by $k_{\min}$.
Consider the optimization problem 
\begin{equation*}
        \min_{x \in E} f(x),
    \end{equation*}
    where $f$ satisfies inequality (\ref{eq:annonying_weaker_smooth_riem}) and is geodesically weak-strongly-convex with parameters $a$ and $\mu$ as in Definition \ref{def:weak-strong-conv_riem}. If $\tilde x$ is produced by one iterate of Riemannian gradient descent (\ref{eq:riemannian_gd}) with $\eta \leq \frac{a}{\zeta \gamma}$, where $\zeta$ is defined as
    \begin{equation*}
        \zeta :=\begin{cases}
    \frac{\sqrt{-k_{\min}} \textnormal{dist}(x,x^*)}{\tanh(\sqrt{-k_{\min}}\textnormal{dist}(x,x^*))} &, k_{\min} < 0 \\
    1 &, k_{\min} \geq 0,
    \end{cases} 
    \end{equation*}
then we have
\begin{equation*}
    \textnormal{dist}^2(\tilde x , x^*) \leq (1-a \mu \eta) \textnormal{dist}^2(x , x^*).
\end{equation*}

\end{prop}

\begin{pf}
    The proof can be found in Appendix B.
\end{pf}

\begin{remark}
As it is evident by previous works in the field [\cite{zhang2016first,alimisis2020continuous}], \emph{convergence is harder in the case of lower curvatures} ($\zeta$ is $1$ if curvatures are nonnegative, but larger than $1$ if curvatures are negative).    
\end{remark}

Before passing to the Riemannian extension of Theorem \ref{thm:main}, we need an auxiliary geometric result similar to Lemma 5 in [\cite{zhang2016first}], which can be found in [\cite{alimisis2020continuous},Corollary 2.1].
\begin{lemma}
\label{le:geometric_bound}
    Let $\Updelta abc$ be a geodesic triangle (i.e. a triangle whose sides are geodesics) in a manifold of sectional curvatures bounded from above by $k_{\max}$. If $k_{\max}>0$, we assume in addition that the lengths of the sides of this triangle are less than $\pi/\sqrt{k_{\max}}$. Then
    \begin{equation*}
     \textnormal{dist}^2(a,c) \geq \delta  \cdot \textnormal{dist}^2(b,c) + 2 \langle \log_{b}(a), \log_{b} (c) \rangle + \textnormal{dist}^2(a,b).
\end{equation*}
where
\begin{equation*}
\delta =\begin{cases}
    \sqrt{k_{\max}} \textnormal{dist}(a,q)\cot(\sqrt{k_{\max}}\textnormal{dist}(a,q)) &, k_{\max} > 0 \\
    1 &, k_{\max} \leq 0,
    \end{cases}
\end{equation*}
with $q$ being some point on the geodesic $bc$.
\end{lemma}
We also need an analogue of Lemma \ref{le:weaker_smoothness}, which can be easily obtained. The reader is referred to the note after Definition 4 in [\cite{zhang2016first}]. 
\begin{lemma}
\label{le:geod_smooth_weaker}
    Let $E$ being an open and geodesically convex subset of $M$. If $f:E \rightarrow \mathbb{R}$ is geodesically $\gamma$-smooth, then it holds
    \begin{equation*}
        \textnormal{dist}^2(x,x^*) \geq \frac{2}{\gamma} (f(x)-f(x^*)).
    \end{equation*}
\end{lemma}
\begin{pf}
    The proof can be found in Appendix B.
\end{pf}

We use these two lemmas to prove a Riemannian analogue of Theorem \ref{thm:main}:
\begin{tcolorbox}
\begin{thm}
\label{thm:main_riemannian}
    Consider the optimization problem
    \begin{equation*}
        \min_{x \in E} f(x),
    \end{equation*}
where $E \subseteq M$ open and geodesically convex with $M$ being a Riemannian manifold of sectional curvatures bounded from above by $k_{\max}$. $f$ is geodesically $\gamma$-smooth and has a unique global minimizer $x^*$. \newline
Assume that a step of Riemannian gradient descent (\ref{eq:riemannian_gd}) starting from any point $x \in E \subseteq M$ satisfies
\begin{equation*}
    \textnormal{dist}^2(\tilde x , x^*) \leq (1-c) \textnormal{dist}^2(x , x^*)
\end{equation*}
for a constant $0<c \leq 1$.
If $k_{\max}>0$, we assume also that $E \subseteq \left \lbrace x \in M / \textnormal{dist}(x,x^*) < \frac{\pi}{4 \sqrt{k_{\max}}} \right \rbrace $ and $\eta \leq \frac{2}{\gamma}$. 

Then $f$ is geodesically weak-strongly-convex in $E$. The parameters of inequality (\ref{eq:weak-quasi-conv_riem}) satisfy
\begin{equation*}
    a := \frac{1}{2 \gamma \eta} \frac{c}{1-\sqrt{\bar \delta c}/2} = \Omega \left(\frac{c}{\gamma \eta} \right), \hspace{2mm}\mu:= \frac{\gamma}{2} = \Omega(\gamma)
\end{equation*}
with
\begin{equation*}
    \bar \delta =\begin{cases}
    \frac{2 \sqrt{k_{\max}} \textnormal{dist}(x,x^*)}{\tan(2 \sqrt{k_{\max}}\textnormal{dist}(x,x^*))} &, k_{\max} > 0 \\
    1 &, k_{\max} \leq 0.
    \end{cases}
\end{equation*}

\end{thm}
\end{tcolorbox}

\begin{pf}
We fix an arbitrary point $x \in E$ and consider $\tilde x$ to be the result of one iterate of Riemannian gradient descent (\ref{eq:riemannian_gd}).

We first bound $\textnormal{dist}^2(\tilde x,x^*)$ using Lemma \ref{le:geometric_bound} in the geodesic triangle $\Updelta x \tilde x x^*$:
    \begin{align*}
        \textnormal{dist}^2(\tilde x,x^*) \geq & \delta \cdot \textnormal{dist}^2(x, \tilde x) +  \textnormal{dist}^2(x,x^*) \\ & - 2 \langle \log_{x}(\tilde x) , \log_x(x^*) \rangle,
    \end{align*}
    where 
\begin{equation*}
\delta =\begin{cases}
\sqrt{k_{\max}} \textnormal{dist}(q,x^*)\cot(\sqrt{k_{\max}}\textnormal{dist}(q,x^*)) &, k_{\max} > 0 \\
1 &, k_{\max} \leq 0,
\end{cases}
\end{equation*}
with $q$ being some point in the geodesic connecting $x$ and $\tilde x$.

This inequality together with the contraction of the assumption gives
\begin{align*}
\label{ineq:1}
    &\delta \cdot \textnormal{dist}^2(x, \tilde x) +  \textnormal{dist}^2(x,x^*) - 2 \langle \log_{x}(\tilde x) , \log_x(x^*) \rangle \leq \\ & (1-c) \textnormal{dist}^2(x,x^*).
\end{align*}
Even when $k_{\max}>0$, $\delta$ can be lower bounded by a positive number $\bar \delta$. 

This is because the function $x \rightarrow x \cot(x)$ is decreasing if $x>0$, thus it suffices to bound $\textnormal{dist}(q,x^*)$ from above.

We do it in the following way:
\begin{equation*}
    \textnormal{dist}(q,x^*) \leq \textnormal{dist}(x,x^*) + \textnormal{dist}(x,q)
\end{equation*}
and \begin{equation*}
    \textnormal{dist}(q,x^*) \leq \textnormal{dist}(\tilde x,x^*) + \textnormal{dist}(\tilde x,q).
\end{equation*}
Thus
\begin{equation*}
    \textnormal{dist}(q,x^*) \leq \frac{\textnormal{dist}(x,x^*)+\textnormal{dist}(\tilde x,x^*)+\textnormal{dist}(x,q)+\textnormal{dist}(\tilde x,q)}{2}.
\end{equation*}
By the assumption, we have that $\textnormal{dist}(\tilde x,x^*) \leq \textnormal{dist}(x,x^*)$ and since $q$ is in the geodesic $x \tilde x$, it holds that $\textnormal{dist}(x,q)+\textnormal{dist}(\tilde x,q) = \textnormal{dist}(x,\tilde x)$.

Since $\tilde x=\exp_x(-\eta \textnormal{grad}f(x))$, we have that $\log_x(\tilde x)=-\eta \textnormal{grad}f(x)$ and $\textnormal{dist}(x,\tilde x) = \eta \|\textnormal{grad}f(x)\|$.
Moreover, by the $\gamma$-smoothness of $f$, we have that $\|\textnormal{grad}f(x)\| \leq \gamma \textnormal{dist}(x,x^*)$.

Using all these facts, we can bound $\textnormal{dist}(q,x^*)$ by
\begin{equation*}
    \textnormal{dist}(q,x^*) \leq \left(1 + \frac{\eta \gamma}{2} \right) \textnormal{dist}(x, x^*) \leq 2 \textnormal{dist}(x,x^*).
\end{equation*}
The last inequality follows from the bound for $\eta$ in the assumption. 

Since $\textnormal{dist}(x,x^*) < \frac{\pi}{4 \sqrt{k_{max}}}$, we have that $2 \textnormal{dist}(x,x^*) < \frac{\pi}{2 \sqrt{k_{max}}}$, thus
\begin{equation*}
    \bar \delta =\begin{cases}
    \frac{2 \sqrt{k_{\max}} \textnormal{dist}(x,x^*)}{\tan(2 \sqrt{k_{\max}}\textnormal{dist}(x,x^*))} &, k_{\max} > 0 \\
    1 &, k_{\max} \leq 0
    \end{cases}
\end{equation*}
is positive and bounds $\delta$ from below.

Using again that $\log_x (\tilde x)=-\eta \textnormal{grad}f(x)$, we can write
\begin{equation}
\label{eq:4}
    2\eta \langle \textnormal{grad}f(x), -\log_x(x^*) \rangle \geq c \cdot \textnormal{dist}^2(x,x^*) + \bar \delta \eta^2 \|\textnormal{grad}f(x)\|^2.
\end{equation}
Next, we use the inequality
\begin{equation*}
    \langle \alpha,\beta \rangle \leq \frac{\rho}{2} \|\alpha\|^2+\frac{1}{2 \rho} \| \beta \|^2
\end{equation*}
for any $\alpha,\beta \in T_x M$ and $\rho>0$ and we obtain
\begin{equation*}
    \frac{\rho}{2} \| \textnormal{grad}f(x) \|^2 \geq \langle \textnormal{grad}f(x), -\log_x(x^*) \rangle-\frac{1}{2 \rho} \textnormal{dist}^2(x,x^*).
\end{equation*}
Multiplying both sides by $\frac{2 \bar \delta \eta^2}{\rho}$, we get
\begin{align*}
    \bar \delta \eta^2 \| \textnormal{grad}f(x)\|^2 & \geq \frac{2 \bar \delta \eta^2}{\rho}\langle \textnormal{grad}f(x), -\log_x(x^*) \rangle \\ &-\frac{\bar \delta \eta^2}{\rho^2} \textnormal{dist}^2(x,x^*).
\end{align*}
Using equation (\ref{eq:4}), we get
\begin{align*}
& 2\eta \langle \textnormal{grad}f(x), -\log_x(x^*) \rangle \geq c \cdot \textnormal{dist}^2(x,x^*)+  \\ & \frac{2 \bar \delta \eta^2}{\rho}\langle \textnormal{grad}f(x), -\log_x(x^*) \rangle-\frac{\bar \delta \eta^2}{\rho^2} \textnormal{dist}^2(x,x^*),
\end{align*}
or equivalently
\begin{align*}
& \left(2\eta - \frac{2 \bar \delta \eta^2}{\rho} \right) \langle \textnormal{grad}f(x),-\log_x(x^*) \rangle \geq \\ & \left( c-\frac{\bar \delta \eta^2}{\rho^2} \right) \textnormal{dist}^2(x,x^*).
\end{align*}
Since the last inequality holds for any $\rho>0$, we can choose $\rho= 2 \frac{\sqrt{\bar \delta} \eta}{\sqrt{c}}$. Then it becomes
\begin{align*}
    & 2 \eta \left( 1-\frac{\sqrt{\bar \delta} \sqrt{c}}{2} \right) \langle \textnormal{grad}f(x),-\log_x(x^*) \rangle \geq \frac{3 c}{4} \textnormal{dist}^2(x,x^*) = \\& \frac{c}{4} \textnormal{dist}^2(x,x^*) + \frac{c}{2} \textnormal{dist}^2(x,x^*).
\end{align*}
By Lemma \ref{le:geod_smooth_weaker}, we have
\begin{equation*}
    \textnormal{dist}^2(x,x^*) \geq \frac{2}{\gamma} (f(x)-f(x^*)),
\end{equation*}
and using that to bound the last term of the previous inequality, we have
\begin{align*}
    & 2 \eta \left( 1-\frac{\sqrt{\bar \delta} \sqrt{c}}{2} \right) \langle \textnormal{grad}f(x),-\log_x(x^*) \rangle \\ &\geq \frac{c}{4} \textnormal{dist}^2(x,x^*) + \frac{c}{\gamma} (f(x)-f(x^*)).
\end{align*}
Rearranging, we get
\begin{align*}
    f(x)-f(x^*) &\leq \underbrace{2 \gamma \eta \frac{ 1-\frac{\sqrt{\bar \delta} \sqrt{c}}{2} }{c 
 }}_{\frac{1}{a}} \langle \textnormal{grad}f(x),-\log_x(x^*) \rangle \\& - \underbrace{\frac{\gamma}{4}}_{\frac{\mu}{2}} \textnormal{dist}^2(x,x^*).
\end{align*}
The last inequality provides the desired result. \qed

\end{pf}

\begin{remark}
Since Theorem \ref{thm:main_riemannian} is in some sense the inverse of Proposition \ref{prop:convergence_riem}, the effect of curvature is the inverse. That is to say, being given a convergence rate in low curvature is a stronger assumption than a convergence rate in higher curvature. The final weak-strong-convexity inequality is therefore weaker in the case of positive curvature ($\delta$ is $1$ for nonpositive curvatures but larger than $1$ for positive ones) and there should be no surprise over that.
\end{remark}

\section{Conclusion}
In this work, we highlighted a deep link between linear convergence of gradient descent and the so-called weak-strong-convexity property. We proved that weak-strong-convexity is sufficient and necessary for obtaining linear convergence guarantees for gradient descent, with the necessary part being the actual novelty. We extended this result to general Riemannian manifolds of bounded sectional curvatures, with the hope that the family of Riemannian gradient descent algorithms parametrised by the chosen Riemannian metric contains many important algorithms. 
If one can show that some arbitrary algorithm for solving a problem with linear convergence can be rewritten as Riemannian gradient descent with respect to some Riemannian metric, then this problem is automatically geodesically weak-strongly-convex with respect to that Riemannian metric. A simple example is Euclidean preconditioned gradient descent (Appendix C). This is not the most fascinating one, but hopefully highlights the capabilities of our theory. We hope to identify more interesting cases in future work.

\begin{ack}
The author would like to thank Bart Vandereycken for proof-reading the main result and offering useful comments on the presenation of the paper. This work is supported by SNSF grant 192363.
\end{ack}

\bibliography{refs}             
                                                   







\newpage

\appendix

\section{Elements of Riemannian geometry}
\label{app:geometry}
In this section we review some basic notions of Riemannian geometry following the exposition of [\cite{alimisis2020continuous}]. This is of course not a complete review and the reader is referred to the classic textbooks for more information. 

\textit{Manifolds:}
A differentiable manifold $M$ is a topological space that is locally Euclidean. This means that for any point $x \in M$, we can find a neighborhood that is diffeomorphic to an open subset of some Euclidean space. This Euclidean space can be proved to have the same dimension, regardless of the chosen point, called the dimension of the manifold.

\textit{Tangent space:}
The tangent space of a differentiable manifold $M$ at a point $x \in M$ is the set of the velocities $\dot \gamma(0)$ of all curves $\gamma:[0,1] \rightarrow M$, such that $\gamma(0)=x$. This space is Euclidean and can be proved to have the same dimension as the one of the manifold.

\textit{Riemannian metrics:}
A Riemannian manifold $(M,g)$ is a differentiable manifold equipped with a Riemannian metric $g_x$, i.e.  an inner product for each tangent space $T_xM$ at $x \in M$. We denote the inner product of $u,v \in T_x M$ with $\langle u,v \rangle_x$ or just $\langle u,v \rangle$ when the tangent space is obvious from context. Similarly we consider the norm at each tangent space as the one induced by the inner product.

\textit{Geodesics:}
Geodesics are curves $\gamma: [0,1] \rightarrow M$ of constant speed and of (locally) minimum length. They can be thought of as the Riemannian generalization of straight lines in Euclidean spaces. A Riemannian manifold where any two points can be connected by some geodesic is called complete. A subset $E \subseteq M$ such that every two points are connected by a unique geodesic is called geodesically convex.

\textit{Exponential map:} Geodesics are used to construct the exponential map $ \textnormal{exp}_x: T_x M \rightarrow M$, defined by $ \textnormal{exp}_x(v)= \gamma(1)$, where $\gamma$ is the unique geodesic such that $\gamma(0)=x$ and $\dot \gamma(0)=v$ (it is a fact that such a geodesic is unique). The exponential map is locally a diffeomorphism.

\textit{Parallel transport:}
Geodesics provide with a natural way to transport vectors from one tangent space to another. This operation is called parallel transport and is usually denoted by $\Gamma_x^y: T_x M \to T_y M$. Details are omitted, but parallel transport is important in order to compare vectors that live in different tangent spaces, as in the case of Definition \ref{def:geod_smoothness}.

\textit{Intrinsic distance:}
Using the notion of geodesics, we can define an intrinsic distance $\textnormal{dist}$ between two points in the Riemannian manifold $M$, as the infimum of the lengths of all geodesics that connect these two points. This notion of distance naturally extends the one from the Euclidean space.

\textit{Riemannian logarithm:} When the exponential map is a diffeomorphism (invertible) its inverse is well-defined. We denote this inverse by $\textnormal{log}_x: M \rightarrow T_x M$. Thus, $\log_x(y)$ denotes the tangent vector at $x$ that one should follow to move from $x$ to $y$ along the geodesic connecting them.

\textit{Curvature:}
In this paper, we make the standard assumption that the input space is not ``infinitely curved". In order to make this statement rigorous, we need the notion of sectional curvatures $K:M \rightarrow \mathbb{R}$, which is a measure of how sharply the manifold is curved (or how ``far" from being flat our manifold is) ``two-dimensionally".

\textit{Riemannian gradient:}
Given a function $f:M \rightarrow \mathbb{R}$, the notions of function differential and Riemannian inner product allow us to define the Riemannian gradient of $f$ at $x \in M$, which is a tangent vector belonging to the tangent space based at $x$, $T_x M$.

\begin{definition}
The Riemannian gradient $\textnormal{\textnormal{gradf}}$ of a (real-valued) function $f:M \rightarrow \mathbb{R}$ at a point $x \in M$, is the tangent vector at $x$, such that $\langle \textnormal{\textnormal{gradf}}(x),u \rangle = df(x)u$~\footnote{$df$ denotes the differential of $f$, i.e. $df (x)[u] = \lim_{t \to 0} \frac{f(c(t)) - f(x)}{t},$ where $c: I \to M$ is a smooth curve such that $c(0) = x$ and $\dot c(0) = u$.}, for any $u \in T_x M$.
\end{definition}

These notions provide all the tools for gradient-based optimization on Riemannian manifolds. They can also provide with a natural extension of the notion of convexity, called geodesic convexity. As full geodesic convexity is not needed in this work, we define only geodesic weak-strong-convexity in the next section.

\newpage

\section{Missing proofs from Section \ref{sec:riemannian}}

\subsection{Proof of Proposition \ref{prop:convergence_riem}}
\begin{pf}
Take an arbitrary $x \in E$ and $\tilde x$ the result of one iterate of Riemannian gradient descent (\ref{eq:riemannian_gd}).

   By Lemma 5 in [\cite{zhang2016first}] combined with Lemma 2 in [\cite{alimisis2020continuous}] (applied to the geodesic triangle $\Updelta x \tilde x x^*$), we have that
   \begin{align*}
       \textnormal{dist}^2(\tilde x,x^*) &\leq \zeta \textnormal{dist}^2(x,\tilde x)+\textnormal{dist}^2(x,x^*) \\&-2 \langle \log_x(\tilde x),\log_x(x^*) \rangle,
   \end{align*}
   where $\zeta$ is 
\begin{equation*}
        \zeta =\begin{cases}
    \frac{\sqrt{-k_{\min}} \textnormal{dist}(x,x^*)}{\tanh(\sqrt{-k_{\min}}\textnormal{dist}(x,x^*))} &, k_{\min} < 0 \\
    1 &, k_{\min} \geq 0.
    \end{cases} 
    \end{equation*}
By noticing that $\log_x(\tilde x)=-\eta \textnormal{grad}f(x)$ by the structure of Riemannian gradient descent, we can rewrite this inequality as
\begin{align*}
   \textnormal{dist}^2(\tilde x,x^*) &\leq \textnormal{dist}^2(x,x^*) -2 \eta \langle \textnormal{grad}f(x), -\log_x(x^*) \rangle \\& +\zeta \eta^2 \| \textnormal{grad}f(x) \|^2.
\end{align*}
By geodesic weak-strong-convexity (Equation (\ref{eq:weak-quasi-conv_riem})), we have that
\begin{align*}
    -2 \eta \langle \textnormal{grad}f(x), -\log_x(x^*) \rangle \leq & -2 \eta a (f(x)-f(x^*)) \\& - a \mu \eta \textnormal{dist}^2(x,x^*).
\end{align*}
Applying inequality (\ref{eq:annonying_weaker_smooth_riem}), we have
\begin{align*}
     -2 \eta \langle \textnormal{grad}f(x), -\log_x(x^*) \rangle \leq & -\frac{\eta a}{\gamma} \|\textnormal{grad}f(x)\|^2 \\& - a \mu \eta \textnormal{dist}^2(x,x^*).
\end{align*}
Plugging that in the bound for $\textnormal{dist}^2(x,x^*)$, we obtain
\begin{align*}
    \textnormal{dist}^2(\tilde x,x^*) & \leq (1-a \mu \eta) \textnormal{dist}^2(x,x^*) \\&+ \left(\zeta \eta^2 -\frac{\eta a}{\gamma} \right) \|\textnormal{grad}f(x)\|^2.
\end{align*}
Since $\eta \leq \frac{a}{\zeta \gamma}$, we have $\zeta \eta^2 -\frac{\eta a}{\gamma} \leq 0$ and the desired result follows.
\end{pf}

\subsection{Proof of Lemma \ref{le:geod_smooth_weaker}}

\begin{pf}
As noted by \cite{zhang2016first} after their Definition 4, the geodesic smoothness inequality (6),
implies that
\begin{equation*}
    f(y) \leq f(x)+\langle \textnormal{grad}f(x) , \log_x(y) \rangle + \frac{\gamma}{2} \textnormal{dist}^2(x,y),
\end{equation*}
for all $x,y \in E$.

For completeness, we provide a detailed proof of this claim:

let $x,y \in E$ and $\beta:[0,1] \rightarrow E$ the unique geodesic (of constant speed) connecting them, such that $\beta(0)=x$ and $\beta(1)=y$. We can write
\begin{align*}
     f(y)-f(x) &=\int_0^1 \frac{d}{dt} f(\beta(t)) dt = \int_0^1 df(\beta(t)) \frac{d}{dt}\beta(t) dt  \\ = & \int_0^1 \langle \textnormal{grad}f(\beta(t)) ,\frac{d}{dt}\beta(t) \rangle dt
\end{align*}
Since $\beta$ is a geodesic of constant speed, it holds
\begin{equation*}
    \frac{d}{dt} \beta(t)=\Gamma_x^{\beta(t)} \log_x(y).
\end{equation*}
Since parallel transport is an isometry and $\Gamma_{p}^q \Gamma_q^p=Id$, we have
\begin{align*}
    & f(y)-f(x)  = \int_0^1 \langle \Gamma_{\beta(t)}^x \textnormal{grad}f(\beta(t)) ,\log_x(y) \rangle dt = \\ &  \langle \textnormal{grad} f(x) , \log_x(y) \rangle  + \int_0^1 \langle \Gamma_{\beta(t)}^x \textnormal{grad}f(\beta(t))-\\ &\textnormal{grad} f(x) ,\log_x(y) \rangle dt.
\end{align*}

The integral can be upper bounded using the Cauchy-Schwarz inequality for the inner product and the assumption (inequality (\ref{eq:geod-smoothness})):
\begin{align*}
    &\int_0^1 \langle \Gamma_{\beta(t)}^x \textnormal{grad}f(\beta(t))-\textnormal{grad} f(x) ,\log_x(y) \rangle dt \leq \\ & \gamma \cdot \textnormal{dist}(x,y) \int_0^1 \textnormal{dist}(\beta(t),x) dt.
\end{align*}
Since $\beta$ is a geodesic of constant speed, we have that
\begin{equation*}
    \textnormal{dist}(\beta(t),x)=t \cdot \textnormal{dist}(x, y),
\end{equation*}
thus
\begin{equation*}
    \int_0^1 \textnormal{dist}(\beta(t),x) dt = \frac{1}{2} \textnormal{dist}(x, y).
\end{equation*}
Putting all together, we indeed get 
\begin{equation*}
    f(y) \leq f(x)+\langle \textnormal{grad}f(x) , \log_x(y) \rangle + \frac{\gamma}{2} \textnormal{dist}^2(x,y).
\end{equation*}

Now we can choose $y=x^*$ and the result follows after noticing that $\textnormal{grad}f(x^*)=0$ and rearranging. \qed
\end{pf}

\newpage
\section{Application to Euclidean preconditioned gradient descent}
\label{app:prec_GD}
We consider the standard Euclidean preconditioned gradient descent algorithm applied to a function $f:\mathbb{R}^n \rightarrow \mathbb{R}$. That is
\begin{tcolorbox}
    \begin{equation}
    \label{eq:prec_gd}
        \tilde x = x-\eta A^{-1} \nabla f(x)
    \end{equation}
\end{tcolorbox}
where $\eta>0$ is a step size and $A$ is a symmetric and positive definite matrix.

This algorithm can be seen as Riemannian gradient descent in $\mathbb{R}^n$ equipped with the Riemannian metric
\begin{equation*}
    \langle u,v \rangle_A := u^T A v.
\end{equation*}
Indeed, since this Riemannian metric does not depend on the specific point of the manifold, the sectional curvatures are $0$ everywhere and the geodesics are still straight lines. This means that
$\exp_x(v)=x+v$ and $\log_x(y)=y-x$.
Moreover, the Riemannian gradient of $f$ is such that
\begin{equation*}
    \langle \textnormal{grad}f(x),v \rangle_A = df(x) v = \langle \nabla f(x),v \rangle,
\end{equation*}
for all $v \in \mathbb{R}^n$.

This can be rewritten to 
\begin{equation*}
    v^T A \textnormal{grad}f(x) = v^T \nabla f(x),
\end{equation*}
and since this holds for all $v$, we have $$\textnormal{grad}f(x)=A^{-1} \nabla f(x).$$

Given the previous discussion, we can apply our Theorem \ref{thm:main_riemannian} to the case of Riemannian gradient descent (\ref{eq:prec_gd}) in the manifold $(\mathbb{R}^n,\langle \cdot,\cdot \rangle_A)$ and obtain the following result:

\begin{tcolorbox}
   \begin{prop}
   \label{prop:prec_gd}
    Consider the optimization problem 
    \begin{equation*}
        \min_{x \in \mathbb{R}^n} f(x),      \end{equation*}
   where $f$ has a unique minimizer $x^*$ and satisfies the inequality
    \begin{equation*}
        \| \nabla f(x)-\nabla f(y) \|_{A^{-1}} \leq \gamma \| x-y \|_A,
    \end{equation*}
    for all $x,y \in \mathbb{R}^n$.

    Then gradient descent preconditioned with a symmetric and positive definite matrix $A$ (\ref{eq:prec_gd}) and with a suitable step size converges with a linear rate
    \begin{equation*}
        \|\tilde x-x^*\|_A^2 \leq (1-c) \|x-x^*\|_A^2,
    \end{equation*}
for some $0<c\leq 1$, if and only if $f$ satisfies the inequality
\begin{align*}
    f(x)-f(x^*) & \leq \frac{1}{a} \langle \textnormal{grad}f(x),x-x^* \rangle \\&-\frac{\mu}{2} (x-x^*)^T A (x-x^*),
\end{align*}
for some $a,\mu>0$ and all $x \in \mathbb{R}^n$.
\end{prop} 

\end{tcolorbox}

\begin{pf}
    As noted before, Euclidean preconditioned gradient descent with $A$ can be written as Riemannian gradient descent with respect to the inner product induced by $A$. We first observe that the inequality
    \begin{equation*}
       \| \nabla f(x)-\nabla f(y) \|_{A^{-1}} \leq \gamma \| x-y \|_A
    \end{equation*}
    is equivalent to geodesic $\gamma$-smoothness with respect to the $A$-inner product, since
    \begin{equation*}
        \| \nabla f(x)-\nabla f(y) \|_{A^{-1}}=\| \textnormal{grad}f(x)-\textnormal{grad}f(y) \|_A.
    \end{equation*}
    
    Thus, applying Proposition \ref{prop:convergence_riem} and restricting the step size to $\eta \leq \frac{a}{\gamma}$, we have that if $f$ is geodesically weak-strongly-convex, then it converges with the mentioned rate, where $c=a \mu \eta$. 
    
  If $f$ converges with the mentioned rate, we can apply Theorem \ref{thm:main_riemannian} and conclude that $f$ must be geodesically weak-strongly-convex with parameters
    \begin{equation*}
    a = \frac{1}{2 \gamma \eta} \frac{c}{1-\sqrt{\bar \delta c}/2} , \hspace{2mm}\mu= \frac{\gamma}{2}.
\end{equation*}
Geodesic weak-strong-convexity in this manifold takes the form
\begin{equation*}
    f(x)-f(x^*) \leq \langle \textnormal{grad}f(x),-\log_x(x^*) \rangle_A -\frac{\mu}{2} \textnormal{dist}^2(x,x^*).
\end{equation*}
For $\langle \textnormal{grad}f(x),-\log_x(x^*) \rangle_A$, we have
\begin{align*}
   &\langle \textnormal{grad}f(x),-\log_x(x^*) \rangle_A =  \langle A^{-1} \nabla f(x), x-x^* \rangle_A = \\ & (x-x^*)^T A A^{-1} \nabla f(x) = \langle \nabla f(x), x-x^* \rangle. 
\end{align*}
For $\textnormal{dist}^2(x,x^*)$, we have
\begin{equation*}
    \textnormal{dist}^2(x,x^*)=\|x-x^*\|_A^2=(x-x^*)^T A (x-x^*).
\end{equation*}
Thus, the desired result follows. \qed
\end{pf}

\begin{remark}
The author is not aware of a prior linear convergence rate with respect to distances of the iterates to the optimum for this algorithm. To the best of our knowledge, the only known result is with respect to function values, see [\cite{uschmajew2022note}, Theorem 2.1].   
\end{remark} 

\begin{remark}
  The conditions of Proposition \ref{prop:prec_gd} are implied by standard Euclidean smoothness and weak-strong-convexity. Indeed, if
$$\| \nabla f(x)-\nabla f(y) \| \leq \gamma \lambda_{\min}(A) \| x-y \|,  $$ then

$$ \frac{1}{\lambda_{min}(A)} \| \nabla f(x)-\nabla f(y) \|^2 \leq \gamma^2 \lambda_{\min}(A) \| x-y \|^2. $$

Now observe that
\begin{align*}
    \| \nabla f(x)-\nabla f(y) \|_{A^{-1}}^2 & \leq \lambda_{\max}(A^{-1}) \| \nabla f(x)-\nabla f(y) \|^2 \\ & = \frac{1}{\lambda_{\min}(A)} \| \nabla f(x) - \nabla f(y) \|^2
\end{align*}
and
\begin{equation*}
\|x-y\|_A^2 \geq \lambda_{\min}(A) \|x-y\|^2.    
\end{equation*}

Thus, it holds
\begin{equation*}
   \| \nabla f(x)-\nabla f(y) \|^2 \leq \gamma^2  \| x-y \|^2.
\end{equation*}
Similarly,
\begin{equation*}
\|x-y\|_A^2 \leq \lambda_{\max}(A) \|x-y\|^2,    
\end{equation*}
thus if
\begin{align*}
    f(x)-f(x^*)  \leq \frac{1}{a} \langle \textnormal{grad}f(x),x-x^* \rangle -\frac{\mu \lambda_{\max}(A)}{2} \|x-x^*\|^2,
    \end{align*}
then
\begin{align*}
    f(x)-f(x^*) \leq \frac{1}{a} \langle \textnormal{grad}f(x),x-x^* \rangle -\frac{\mu}{2} (x-x^*)^T A (x-x^*).
\end{align*}
  
\end{remark}

\end{document}